\documentclass[10pt]{article}
\usepackage{amsmath, amsthm}
\usepackage{amsfonts}
\usepackage{indentfirst}
\newtheoremstyle{theorem}
{10pt} 
{10pt} 
{\sl} 
{\parindent} 
{\bf} 
{. } 
{ } 
{} 
\theoremstyle{theorem}

\newtheorem{propos}{Proposition}
\newtheoremstyle{defi}
{10pt} 
{10pt} 
{\rm} 
{\parindent} 
{\bf} 
{. } 
{ } 
{} 
\theoremstyle{defi}

\begin{document}

\title{Bifurcation and numerical study in an EHD convection problem}
\author{Ioana Dragomirescu\\
Dept. of Math.,
 Univ. "Politehnica" of Timisoara,\\
 Piata Victoriei, No.2, 300006,
  Timisoara, Romania\\
 i.dragomirescu@gmail.com\\
}
\date{}
\maketitle

\begin{abstract}
The linear eigenvalue problem governing the stability of the
mechanical equilibrium of the fluid in a electrohydrodynamic (EHD)
convection problem is investigated. The analytical study is one of
bifurcation. This allows us to regain the expression of the
neutral surface in the classical case. The method used in the
numerical study is a Galerkin type spectral method based on
polynomials and it provides good results.
\end{abstract}

\section{The physical problem}

The investigated physical model is one of the two EHD models of
Roberts \cite{R1}, based on the Gross'experiments \cite{Gr} which
are concerned with a layer of insulating oil confined between two
horizontal conducting planes and heated from above and cooled from
below. The experimental investigations showed that the presence of
a vertical electric field of sufficient strength across the layer,
lead to a tesselated pattern of motions, in a manner similar to
that of the classical B\'{e}nard convection \cite{Str}. Gross
\cite{Gr} suggested that this phenomenon may be due to the
variation of the dielectric constant of the fluid with the
temperature.
\par In the first model investigated by Roberts \cite{R1} the dielectric
constant is allowed to vary with the temperature. The homogeneous
insulating fluid is assumed to be situated in a layer of depth $d$
(the fluid occupies the region between the planes $z=\pm 0.5d$,
which are maintained at uniform but different temperatures), with
vertical, parallel applied gradients of temperature and
electrostatic potential. The uniform electric field is applied in
the $z$ direction.
\par The equations governing the EHD convection upon normal mode
representations are \cite{Str}
\begin{equation}
\label{eq: model1_roberts}
(D^{2}-a^{2}-s)(D^{2}-a^{2}-Prs)(D^{2}-a^{2})^{2}F=La^{4}F-R_{a}a^{2}(D^{2}-a^{2})F
\end{equation}
with the boundary conditions on $F$
\begin{equation}
\label{eq:bc_model1}
\begin{array}{l}
F=D^{2}F=D(D^{2}-a^{2}-Pr s)F=0\\
\\
((D^{2}-a^{2})(D^{2}-a^{2}-s)(D^{2}-a^{2}-Pr s)+R_{a}a^{2})(DF\pm
kaF)=0
\end{array} \textrm{ at }z=\pm 0.5.
\end{equation}
\par Here, the unknown function $F$ is the amplitude of the temperature perturbation $\Theta$, i.
e. $\Theta=F(z)e^{(i(lx+my)+st)}$, $Pr$ is the Prandtl number,
$k=\dfrac{\epsilon_{m}}{\widehat{\epsilon}}$, with $\epsilon_{m}$
the value of the dielectric field at the temperature $T_{m}=0$ and
$\widehat{\epsilon}$ the electric constant of the solid in
$z>\dfrac{1}{2}$, $a$ is the wavenumber, $a^{2}=l^{2}+m^{2}$,
$R_{a}$ is the Rayleigh number, $L$ is a parameter measuring the
potential difference between the planes.
\par Roberts \cite{R1} investigated only the stationary case, i.e. $s=0$,
so the eigenvalue problem consists from an eight-order
differential equation
\begin{equation}
\label{eq:model1}
(D^{2}-a^{2})^{4}F-La^{4}F+R_{a}a^{2}(D^{2}-a^{2})F=0
\end{equation}
and the boundary conditions
\begin{equation}
\label{eq:bc_m1}
F=D^{2}F=D(D^{2}-a^{2})F=((D^{2}-a^{2})^{3}+R_{a}a^{2})(DF\pm
kaF)=0 \textrm{ at }z=\pm 0.5.
\end{equation}
He found that when the smallest Rayleigh number,
$R_{a_{\min}}=\min _{a}R_{a}(a)$, varies from $-1000$ to
$1707.762$, $L$ decreases from $3370.077$ to 0.
\par The second model \cite{R1} was also been investigated by
Turnbull \cite{T1}, \cite{T2}. In this case, the variation of the
dielectric constant is not important but the fluid is weakly
conducting and its conductivity varies with temperature. \par The
eigenvalue equation has the form \cite{Str}
\begin{equation}
\label{eq:model2} (D^{2}-a^{2})^{3}F+R_{a}a^{2}F+Ma^{2}DF=0
\end{equation}
with $M$ a dimensionless parameter measuring the variation of the
electrical conductivity with temperature. The boundary conditions,
written for the case of rigid boundaries at constant temperatures,
read
\begin{equation}
\label{eq:bc_m2} F=D^{2}F=D(D^{2}-a^{2})F=0 \textrm { at } z=\pm
0.5
\end{equation}
Roberts \cite{R1} found that when $M$ is increasing from 0 to
1000, $R_{a_{\min}}$ is increasing from 1707.062 to 2065.034.
\par Straughan \cite{Str} also investigated  these EHD convection
problems, developing a fully nonlinear energy stability analysis
for non-isothermal convection problems in a dielectric fluid.
\section{The bifurcation analysis}
The linear stability of the motion or the equilibrium of a fluid
in many problem from hydrodynamic, electrohydrodynamic or
hydromagnetic stability theory is governed by a linear
higher-order ordinary differential equation  with constant
coefficients and homogeneous boundary conditions. The exact
solution of such equations or, for the case of eigenvalue
problems, the exact eigenvalue  is most of the times impossible to
find. That is why, numerical methods, usually implying an infinite
number of terms, leading however to an approximative solution by
some specific truncations to a finite number of terms, are used.
However, the theoretical methods can impose restrictions with
regards to the numerical results.
 \par For the considered problem let us
 introduce the direct method \cite{Ge1} which consists in the determination of the
 eigenfunctions and their introduction into the boundary
 conditions.
 \par The characteristic equation associated to (\ref{eq:model2})
 is
 \begin{equation}
 \label{eq:ec_car}
 (\lambda^{2}-a^{2})^{3}+Ma^{2}\lambda+R_{a}a^{2}=0
 \end{equation}
 When the characteristic equation has multiple roots the straightforward application
of numerical method can lead to false secular points. That is why,
these cases must be investigated separately.
\begin{propos}
For $M=0$  the only secular points are those situated on
$$NS_{n}:R_{a}=\dfrac{((2n-1)^{2}\pi^{2}+a^{2})^{3}}{a^{2}}, \ \forall n\in \mathbb{N}.$$
\end{propos}
\begin{proof}
\par For $M=0$, the characteristic equation (\ref{eq:ec_car})
reduces to
\begin{equation}
\label{eq:ec_car_red} (\lambda^{2}-a^{2})^{3}+a_{2}=0, \textrm {
with } a_{2}=R_{a}a^{2}.
\end{equation}
In this classical case the roots of (\ref{eq:ec_car_red}) have the
form
$$\begin{array}{l}
\lambda_{1,2}=\sqrt{a^{2}+\sqrt[3]{-a_{2}}\epsilon_{1,2}}, \
\lambda_{3}=\sqrt{a^{2}+\sqrt[3]{-a_{2}}},\\
\\
\lambda_{4}=-\lambda_{1}, \lambda_{5}=-\lambda_{2},
\lambda_{6}=-\lambda_{3}
\end{array}$$
so the general solution of (\ref{eq:model2}) has the form
$$F=\sum\limits_{i=1}^{3}A_{i}\cosh(\lambda_{i}z)+B_{i}\sinh(\lambda_{i}z).$$
Replacing the solution $F$ into the boundary conditions
(\ref{eq:bc_m2}) we get the secular equation
\begin{equation}
\label{eq:ec_sec} \Delta= \begin{tabular}{|cccccc|}
$0$&$0$&$0$&$m_{1}$&$m_{2}$&$m_{3}$\\
$1$&$1$&$1$&$0$&$0$&$0$\\
$0$&$0$&$0$&$\lambda_{1}^{2}m_{1}$&$\lambda_{2}^{2}m_{2}$&$\lambda_{3}^{2}m_{3}$\\
$\lambda_{1}^{2}$&$\lambda_{2}^{2}$&$\lambda_{3}^{2}$&$0$&$0$&$0$\\
$-\lambda_{1}\mu_{1}m_{1}$&$-\lambda_{2}\mu_{2}m_{2}$&$-\lambda_{3}\mu_{3}m_{3}$&$0$&$0$&$0$\\
$0$&$0$&$0$&$\lambda_{1}\mu_{1}$&$\lambda_{2}\mu_{2}$&$\lambda_{3}\mu_{3}$\\
\end{tabular}=0
\end{equation}
with $m_{i}=\tanh(\lambda_{i}/2)$,
$\mu_{i}=\lambda_{i}^{2}-a^{2}$, $i=1,2,3$. \par When
$\cosh(\lambda_{i}/2)\neq 0$, $i=1,2,3$, we can rewrite the
secular equation as $\Delta=\Delta_{1}\cdot \Delta_{2}$ with
$$\Delta_{1}=\lambda_{1}\mu_{1}m_{1}(\lambda_{2}^{2}-\lambda_{3}^{2})+\lambda_{2}\mu_{2}m_{2}(\lambda_{3}^{2}-\lambda_{1}^{2})+
\lambda_{3}\mu_{3}m_{3}(\lambda_{1}^{2}-\lambda_{2}^{2})$$ and
$$\Delta_{2}=\lambda_{1}^{2}m_{1}(\lambda_{3}\mu_{3}m_{2}-\lambda_{2}\mu_{2}m_{2})+\lambda_{2}^{2}m_{2}(\lambda_{1}\mu_{1}m_{3}-\lambda_{3}\mu_{3}m_{1})+
\lambda_{3}^{2}m_{3}(\lambda_{2}\mu_{2}m_{1}-\lambda_{1}\mu_{1}m_{2}).$$
For $a>0$, the equations $\Delta_{1}=0$ and $\Delta_{2}=0$ have
only null solutions $R=0$, so no secular points exists on these
surfaces.
\par The condition $\cosh(\lambda_{i}/2)\neq 0$, $i=1,2,3$ is not
fulfilled  only for $i=3$, i.e.
$\cosh(\lambda_{3}/2)=0\Leftrightarrow
\cos(\lambda_{3}/2)=0\Leftrightarrow
\lambda_{3}^{2}=-(2n-1)^{2}\pi^{2}$, which implies that the
secular curve is $NS_{n}$. And, indeed, the critical values of the
Rayleigh number $R_{a}$
 belong to $NS_{1}$ only, identical to the classical one from Chandrasekhar
 \cite{Ch1}.
 \end{proof}
 \par The general form of the solution of the two-point problem for the governing differential equation
is written in terms of the roots of the characteristic equation
associated with the differential equation. In addition, this form
depends on the multiplicity of the characteristic roots.
Introducing the general solution into the boundary conditions the
secular equation is obtained and it depends on the multiplicity of
the characteristic roots. As a consequence, the secular equation
has different forms in different regions of the parameter space.
Each eigenvalue is a solution of the obtained secular equation, so
the eigenvalue depends on all other physical parameters. The
neutral manifolds (the most convenient manifolds from the physical
point of view), generated by the secular equation separate the
domain of stability from the domain of instability.
\par Let us consider the general case when the roots of the characteristic equation $\lambda_{1}$,
 $\lambda_{2}$, ..., $\lambda_{6}$ are distinct. Then the general solution of (\ref{eq:model2})
 has the form $
 F(z)=\sum\limits_{i=1}^{6}A_{i}e^{\lambda_{i}z}.
$ Introducing it into the boundary conditions (\ref{eq:bc_m2}) we
obtain the secular equation \cite{Drag1}
\begin{equation}
\label{eq:ec_sec_generala} \Delta(a, M,R_{a})=0,
\end{equation}
where $\Delta$ is a determinant. Its  $i$-th column has the same
form in $\lambda_{i}$ as any other $j$-th column in $\lambda_{j}$.
If $\lambda_{i}=\lambda_{j}$, then the $i$-th and the $j$-th
columns in $\Delta$ are identical. Therefore $\Delta\equiv 0$. In
fact, in this situation, (\ref{eq:ec_sec_generala}) is not
entitled to serve as a secular equation and the direct numerical
computations will be invalid. When $M\neq 0$, some particular
cases interesting from the bifurcation point of
 view,
 arise due to the existence of bifurcation sets of the characteristic manifold.\par

 Let us consider  the surface
 $S_{0}$ defined by the points $(a,M,R_{a})=(a,M,a^{4})$. In this case
 we have the following result
\begin{propos}
Let us define the surfaces
$$
S_{i}:
M=\dfrac{(33\mp3\sqrt{21})\sqrt{90\pm10\sqrt{21}}a^{3}}{250},\
i=1,2.
$$
The surface $S_{0}\cap S_{i}$, $i=1,2$ is a bifurcation set of the
characteristic manifold defined by (\ref{eq:ec_car}). The points
on $S_{0}\cap S_{i}$, $i=1,2$ are not secular.
\end{propos}
\begin{proof}
If $(a,M,R_{a})\in S_{0}$ then $R_{a}=a^{4}$ and one of the roots
of the characteristic equation is , for instance,
$\lambda_{1}=0$. Assuming that $M\neq 0$ and $a>0$, $\lambda_{1}$
is not a double root of (\ref{eq:ec_car}). The search of multiple
roots reduces then to the equation
\begin{equation}
\label{eq:ec_car_0}
\lambda^{5}-3\lambda^{4}a^{2}+3\lambda^{2}a^{4}+Ma^{2}=0.
\end{equation}
No multiple roots of algebraic multiplicity order greater than 2
exists. The double roots of (\ref{eq:ec_car_0}) must also be roots
of its derivative $5z^{4}-9z^{2}a^{2}+3a^{4}=0$. In these
conditions the possible  double roots are $\lambda_{2,3}=\lambda=-
\dfrac{1}{10}\sqrt{90\pm 10\sqrt{21}}a$ only for $(a,M, R_{a})\in
S_{i}$, $i=1,2$, i.e. the surfaces $S_{0}\cap S_{i}$, $i=1,2$ are
bifurcation sets for the characteristic equation
(\ref{eq:ec_car}). In the case of multiple roots, the general form
of the solution of (\ref{eq:model2}) is
$F(z)=\sum\limits_{i=1}^{n}P_{i}(z)e^{\lambda_{i}z}$, where
$P_{i}$ is an algebraic polynomial of $m_{i}-1$ degree, $m_{i}$
being the algebraic multiplicity of $\lambda_{i}$, in our case
$F(z)=A+(B+Cz)e^{\lambda}z+\sum\limits_{i=4}^{6}A_{i}e^{\lambda_{i}z}$.\par
Formally, the secular equation is deduced from
(\ref{eq:ec_sec_generala}), by writing the column $i$ for
$\lambda_{i}$ while the columns $i+1$, $i+2$,..., $i+m_{i}-1$ are
obtained by differentiating $l$, $l=1,2,...,m_{i}-1$ times the
$i+l$-th column of (\ref{eq:ec_sec_generala}) with respect to
$\lambda_{i+l}$ and then replacing $\lambda_{i+l}$ by
$\lambda_{i}$.
\par However, the numerical evaluations show that now secular
points exists on these surfaces.
\par
\end{proof}
\section{Spectral methods based study}
\par The second part of our study regards the numerical treatment of the two-point problem (\ref{eq:model2}) - (\ref{eq:bc_m2}).
\par Large classes of eigenvalue problems can be solved numerically
using spectral methods, where, typically, the various unknown
fields are expanded upon sets of orthogonal polynomials or
functions. The convergence of such methods is in most cases easy
to assure and they are efficient, accurate and fast. Our numerical
study is performed using a weighted residual (Galerkin type)
spectral method.
\par  Introducing the
new function $U=(D^{2}-a^{2})F$, the generalized eigenvalue
problem
\begin{equation}
\label{eq:gen_eig} \left\{
\begin{array}{l}
(D^{2}-a^{2})^{2}U=-R_{a}a^{2}F-Ma^{2}DF,\\
(D^{2}-a^{2})F=U,\\
F=U=DU=0 \textrm { at }z=\pm 0.5.
\end{array}
\right.
\end{equation}
is obtained. Following \cite{HSTR}, we consider the orthogonal
sets of functions
$$
\begin{array}{l}
\{\phi_{i}\}_{i=1,2,...,N}:
\phi_{i}(z)=\int_{-0.5}^{z}L_{i}^{*}(t)dt, \textrm { verifying } \phi_{i}(\pm 0.5)=0,\\
\\
\{\beta_{i}\}_{i=1,2,...,N}:
\beta_{i}(z)=\int_{-0.5}^{z}\int_{-0.5}^{s}L_{i}^{*}(t)dtds,
\textrm { verifying } \beta_{i}(\pm 0.5)=D\beta_{i}(\pm 0.5)=0,
\end{array}
$$
with $L_{i}^{*}=L_{i}(2x)$ the shifted Legendre polynomials on
$(-0.5,0.5)$ and $L_{i}$ the Legendre polynomials on $(-1,1)$. The
unknown functions from (\ref{eq:gen_eig}), $U,F$, are written as
truncated series of functions $\beta_{i}$, irrespective
$\phi_{i}$, i.e.
$$U=\sum\limits_{i=1}^{N}U_{i}\beta_{i}(z), \ F=\sum\limits_{i=1}^{N}F_{i}\phi_{i}(z).$$
The boundary conditions on $U$ and $F$ are automatically
satisfied. Replacing these expressions in (\ref{eq:gen_eig}),
imposing the condition of orthogonality on the vector
$(\beta_{k},\phi_{k})^{T}$, $k=1,2,...,N$, we get an algebraic
system in the unknown, but not all null, coefficients $U_{i}$,
$F_{i}$. The secular equation, written as the determinant of the
obtained algebraic system, gives us the values of the Rayleigh
number as a function of the other physical parameters. The
smallest values of the Rayleigh number for various values of the
parameters $a$ and $M$ form the neutral surface that separates the
domain of stability from the instability domain. All the
expression of the scalar products resulting in the algebraic
system are given in \cite{Drag2} for the case of shifted Legendre
polynomials on $(0,1)$, but they are easy to adjust to the
interval $(-0.5..0.5)$. The specific choice of basis functions led
to sparse matrices, with banded sub-matrices of dimension $N\times
N$.
\par The numerical evaluations of the critical Rayleigh numbers were
obtained for a small number of terms $N$ ($N=6$) in the truncated
series confirming the well-known accuracy of spectral methods. We
obtained that critical values of $R_{a}$ are increasing from
$1734.120$ to $2082.808$ when $M$ is increasing from $0$ to
$1000$, similar to the ones of Roberts \cite{R1}.
\par The unknown vector fields from (\ref{eq:gen_eig}) can also be
expanded upon complete sequences of functions in $L^{2}(-0.5,0.5)$
defined by using Chebyshev polynomials that satisfy the boundary
conditions of the problem. Keeping the above notations, the
functions $\phi_{i}$, $i=1,2,...,N$ are defined by
$\phi_{i}(z)=T^{*}_{i}(z)-T^{*}_{i+2}(z)$ and $\beta_{i}$,
$i=1,2,...,N$ by
$\beta_{i}(z)=T^{*}_{i}-\dfrac{2(i+2)}{i+3}T^{*}_{i+2}+\dfrac{i+1}{i+3}T^{*}_{i+4}$
\cite{Sh2} with $T^{*}_{i}$, $i=1,2,...,N$, the shifted Chebyshev
polynomials on $(-0.5,0.5)$ defined in a similar manner as the
shifted Legendre polynomials. All the evaluations of the scalar
products were based on the orthogonality relation
 \begin{equation}
 \label{eq:CG_ort}
 \int_{-0.5}^{0.5}T_{n}^{*}(z)T_{m}^{*}(z)w^{*}(z)dz=\left\{
 \begin{array}{l}
 \dfrac{\pi}{2}c_{n}\delta_{nm}, \textrm{ if } i=j,\\
0, \textrm{ if } i\neq j, \end{array}
 \right.
 \end{equation}
 with respect to the weight function $w^{*}(z)=\dfrac{1}{\sqrt{1/4-z^{2})}}$.
\par The numerical results where
obtained for a larger number of terms in the expansion sets
($N=11$) and they show that the shifted Legendre based method is
more effective in this case. We can mention that the Chebyshev
polynomials are considered suitable more likely for the tau method
or the collocation type methods. Some numerical evaluations of
$R_{a}$ as a function of $a$ and $M$ are given in Table 1.
\par Other sets of complete orthogonal functions based on
Chebyshev polynomials and satisfying various boundary conditions
can be found in \cite{Gh}, \cite{Mas1}.
\begin{table}[tbp] \centering%
\begin{tabular}{|c|c|c|c|}
\hline  $a$ & $M$ & $R_{a}-SLP$&$R_{a}-SCP$ \\
\hline  $3.117$ & $0$ & $1734.120$ &$1775.955$\\
\hline $3.117$&$10$ & $1734.154$ & $1775.987$ \\
\hline $3.117$ & $1000$ & $2082.802$&$2100.935$\\
 \hline $1.5$ & $0$ & $3116.286$ & $31199.286$ \\
\hline $1.5$ & $5$ & $3116.381$ & $3199.289$ \\
\hline $10$ & $0$ & $11409.157$ & $14909.559$ \\
\hline $10$ & $100$ & $14414.05$ &
$14419.963$ \\
 \hline $10$ & $500$ & $14531.694$ & $14994.747$ \\
\hline $20$ & $0$ & $166779.036$ & $182878.881$ \\ \hline
\end{tabular}%
\caption{Numerical values for the Rayleigh number for various
values of the parameters $a$, $M$ obtained by spectral methods
based on shifted Legendre (SLP) and shifted Chebyshev (SCP)
polynomials. \label{V4}}%
\end{table}%

 \section{Conclusions}
 In this paper  we performed a bifurcation analysis and a numerical treatment for an electrohydrodynamic convection problem.
 When eigenvalue problems from linear stability theory are
 investigated only numerically spurious eigenvalue can be
 encountered, especially when bifurcation sets of the characteristic manifold occur. In order to detect the false secular points a
 bifurcation study of the problem becomes necessary. An example of this type of
 problems was investigated in \cite{Ge2}, e.g. for an electrohydrodynamic
 convection problem in the case of free-free boundaries the numerical methods led to the existence of false secular
 points.
 \par The numerical study was performed  here using a Galerkin
 type spectral method which implied that the boundary conditions
 are satisfied by the orthogonal sets of expansion functions.
 However, when this condition is not fulfilled, the tau method or
 the collocation method can also be applied.  All these methods
 are widely used in the numerical investigation of eigenvalue problems
 governing the linear stability of motions or equilibrium of
 fluid in convection problems.  From the physical point of view, the evaluations of the Rayleigh number $R_{a}$ showed
 an enlargement of the stability domain when the parameter $M$ is
 increasing, the dependence of $R_{a}$ of $M$ is not however exponential. These evaluations were easy to compute and proved to be similar to the ones
 existing in the literature.

\end{document}